\documentclass[12pt]{article}
\usepackage{oldlfont}
\usepackage{enumerate}
\usepackage{amsmath}
\usepackage{amssymb}
\usepackage{amsthm}
\usepackage{mathrsfs}
\usepackage{amscd}
\usepackage{exscale}
\usepackage{latexsym}
\usepackage[all]{xy}
\usepackage{hyperref}
\usepackage{fancyhdr}
\usepackage{layout}
\usepackage{color}
\usepackage{graphicx}
\DeclareFontFamily{U}{wncy}{}
\DeclareFontShape{U}{wncy}{m}{n}{<->wncyr10}{}
\DeclareSymbolFont{mcy}{U}{wncy}{m}{n}
\DeclareMathSymbol{\Sh}{\mathord}{mcy}{"58} 

\begin{document}

\baselineskip=17pt

\pagestyle{headings}

\numberwithin{equation}{section}

\makeatletter                                                           

\def\section{\@startsection {section}{1}{\z@}{-5.5ex plus -.5ex         
minus -.2ex}{1ex plus .2ex}{\large \bf}}                                 


\pagestyle{fancy}
\renewcommand{\sectionmark}[1]{\markboth{ #1}{ #1}}
\renewcommand{\subsectionmark}[1]{\markright{ #1}}
\fancyhf{} 
\fancyhead[LE,RO]{\slshape\thepage}
\fancyhead[LO]{\slshape\rightmark}
\fancyhead[RE]{\slshape\leftmark}

\addtolength{\headheight}{0.5pt} 
\renewcommand{\headrulewidth}{0pt} 

\newtheorem{thm}{Theorem}
\newtheorem{mainthm}[thm]{Main Theorem}

\newcommand{\ZZ}{{\mathbb Z}}
\newcommand{\GG}{{\mathbb G}}
\newcommand{\Z}{{\mathbb Z}}
\newcommand{\RR}{{\mathbb R}}
\newcommand{\NN}{{\mathbb N}}
\newcommand{\GF}{{\rm GF}}
\newcommand{\QQ}{{\mathbb Q}}
\newcommand{\CC}{{\mathbb C}}
\newcommand{\FF}{{\mathbb F}}

\newtheorem{lem}[thm]{Lemma}
\newtheorem{cor}[thm]{Corollary}
\newtheorem{pro}[thm]{Proposition}
\newtheorem{proprieta}[thm]{Property}
\newcommand{\pf}{\noindent \textbf{Proof.} \ }
\newcommand{\eop}{$_{\Box}$  \relax}
\newtheorem{obs}[thm]{Remark}
\newtheorem{num}{equation}{}

\theoremstyle{definition}
\newtheorem{rem}[thm]{Remark}
\newtheorem*{D}{Definition}

\newcommand{\nsplit}{\cdot}
\newcommand{\G}{{\mathfrak g}}
\newcommand{\GL}{{\rm GL}}
\newcommand{\SL}{{\rm SL}}
\newcommand{\SP}{{\rm Sp}}
\newcommand{\LL}{{\rm L}}
\newcommand{\Ker}{{\rm Ker}}
\newcommand{\la}{\langle}
\newcommand{\ra}{\rangle}
\newcommand{\PSp}{{\rm PSp}}
\newcommand{\U}{{\rm U}}
\newcommand{\GU}{{\rm GU}}
\newcommand{\Aut}{{\rm Aut}}
\newcommand{\Alt}{{\rm Alt}}
\newcommand{\Sym}{{\rm Sym}}

\newcommand{\isom}{{\cong}}
\newcommand{\z}{{\zeta}}
\newcommand{\Gal}{{\rm Gal}}

\newcommand{\F}{{\mathbb F}}
\renewcommand{\O}{{\cal O}}
\newcommand{\Q}{{\mathbb Q}}
\newcommand{\R}{{\mathbb R}}
\newcommand{\N}{{\mathbb N}}
\newcommand{\A}{{\mathcal{A}}}
\newcommand{\E}{{\mathcal{E}}}
\newcommand{\J}{{\mathcal{J}}}


\newcommand{\DIM}{{\smallskip\noindent{\bf Proof.}\quad}}
\newcommand{\CVD}{\begin{flushright}$\square$\end{flushright}
\vskip 0.2cm\goodbreak}


\vskip 0.5cm

\title{Counterexamples to the local-global divisibility over elliptic curves}
\author{Gabriele Ranieri\footnote{Instituto de Matem\'aticas, Pontificia Universidad Cat\'olica de Valpara\'iso, supported by Fondecyt Regular Project N. 1140946}}
\date{  }
\maketitle

\vskip 1.5cm

\begin{abstract}
Let $p \geq 5$ be a prime number. We find all the possible subgroups $G$ of ${\rm GL}_2 ( \Z / p \Z )$ such that there exists a number field $k$ and an elliptic curve $\E$ defined over $k$ such that the $\Gal ( k ( \E[p] ) / k )$-module $\E[p]$ is isomorphic to the $G$-module $( \Z / p \Z )^2$ and there exists $n \in \N$ such that the local-global divisibility by $p^n$ does not hold over $\E ( k )$.  
\end{abstract}

\section{Introduction}

Let $k$ be a number field and let ${\mathcal{A}}$ be a commutative algebraic group defined over $ k $.
Several papers have been written on  the following classical question, known as the \emph{Local-Global Divisibility Problem}.

\par\bigskip\noindent  P{\small ROBLEM}: \emph{Let $P \in {\mathcal{A}}( k )$. Assume that for all but finitely many valuations $v$ of $k$, there exists $D_v \in {\mathcal{A}}( k_v )$ such that $P = qD_v$, where $q$ is a positive integer. Is it possible to conclude that there exists $D\in {\mathcal{A}}( k )$ such that $P=qD$?}

\par\bigskip\noindent  By  B\'{e}zout's identity, to get answers for a general integer it is sufficient to solve it for powers $p^n$ of a prime. In the classical case of ${\mathcal{A}}={\mathbb{G}}_m$, the answer is positive for $ p $ odd, and negative for instance for $q=8$ (and $P=16$) (see for example \cite{AT}, \cite{Tro}).

\bigskip  For general commutative algebraic groups, Dvornicich and Zannier gave a cohomological interpretation of the problem (see \cite{DZ1} and \cite{DZ3}) that we shall explain.
Let $ \Gamma $ be a group and let $ M $ be a $\Gamma$-module.
We say that a cocycle $Z \colon \Gamma \rightarrow M$ satisfies the local conditions if for every $\gamma \in \Gamma$ there exists $ m_\gamma \in M $ such that 
$Z_\gamma = \gamma ( m_\gamma ) - m_\gamma$.
The set of the class of cocycles in $H^1 ( \Gamma, M )$ that satisfy the local conditions is a subgroup of $H^1 ( \Gamma, M )$.
We call it the first local cohomology group $H^1_{{\rm loc}} ( \Gamma, M )$. 
Dvornicich and Zannier \cite[Proposion 2.1]{DZ1} proved the following result.

\begin{pro}\label{pro1}
Let $ p $ be a prime number, let $n$ be a positive integer, let $k$ be a number field and let $\A$ be a commutative algebraic group defined over $k$.
If $H^1_{{\rm loc}} ( \Gal ( k ( \A[p^n] ) / k ) , \A[p^n] ) = 0$, then the local-global divisibility by $p^n$ over $\A ( k )$ holds.
\end{pro}

The converse of Proposition \ref{pro1} is not true, but in the case when the group $H^1_{{\rm loc}} ( \Gal ( k ( \A[p^n] ) / k ) , \A[p^n] )$ is not trivial we can find an extension $ L $ of $k$ such that $L \cap k ( \A[p^n] ) = k$, and the local-global divisibility by $p^n$ over $\A ( L )$ does not hold (see \cite[Theorem 3]{DZ3} for the details).

Many mathematicians got criterions for the validity of the local-global divisibility principle for many families of commutative algebraic groups, as algebraic tori (\cite{DZ1} and \cite{Ill}), elliptic curves (\cite{Cre1}, \cite{Cre2}, \cite{DZ1}, \cite{DZ2}, \cite{DZ3}, \cite{GR1}, \cite{Pal1}, \cite{Pal2}, \cite{PRV1}, \cite{PRV2}), and very recently polarized abelian surfaces (\cite{GR2}) and ${\rm GL}_2$-type varieties (\cite{GR3}).

In this paper we are interested in the family of the elliptic curves.
Let $p$ be a prime number, let $ k $ be a number field and let $\E$ be an elliptic curve defined over $k$. 
Dvornicich and Zannier \cite[Theorem 1]{DZ3} found a very interesting criterion for the validity of the local-global divisibility by a power of $p$ over $\E ( k )$, in the case when $k \cap \Q ( \zeta_p ) = \Q$.

In a joint work with Paladino and Viada (see \cite{PRV1}, and Section \ref{sec11}) we refined this criterion, by proving that if $ k $ does not contain $\Q ( \zeta_p + \overline{\zeta_p} )$ and $\E ( k )$ does not admit a point of order $p$, then for every positive integer $n$, the local-global divisibility by $p^n$ holds over $\E ( k )$.
In another joint work with Paladino and Viada \cite{PRV2} we improved our previous criterion and the new criterion allowed us to show that if $k = \Q$ and $p \geq 5$, for every positive integer $n$ the local-global divisibility by $p^n$ holds for $\E ( \Q )$.     

Very recently, Lawson and Wutrich \cite{LW} found a very strong criterion for the triviality of $H^1 ( \Gal ( k ( \E[p^n] ) / k ) , \E[p^n] )$ (then for the validity of the local-global principle by $p^n$ over $\E ( k )$, see Proposition \ref{pro1}), but still in the case when $k \cap \Q ( \zeta_p ) = \Q$.

Finally, Dvornicich and Zannier \cite{DZ2} and Paladino (\cite{Pal1}) studied the case when $p = 2$ and Paladino (\cite{Pal2}) and Creutz (\cite{Cre1}) studied the case when $p = 3$.

Then, we have a fairly good understanding of the local-global divisibility by a power of $p$ over $\E ( k )$ either when $p \in \{ 2, 3 \}$ or $k$ does not contain $\Q ( \zeta_p + \overline{\zeta_p} )$ and $\E ( k )$ does not admit a point of order $p$.  
In this paper we study the cases not treated, yet. 
We get the following results:

\begin{thm}\label{teo1}
Let $p \geq 5$ be a prime number, let $ k $ be a number field and let $\E$ be an elliptic curve defined over $k$. Suppose that there exists a positive integer $n$ such that the local-global divisibility by $p^n$ does not hold over $\E ( k )$. Let $G_1$ be $\Gal ( k ( \E[p] ) / k )$.
Then one of the following holds:
\begin{enumerate}
\item $G_1$ is cyclic of order dividing $p-1$ generated by an element with an eigenvalue equal to $1$;
\item $p \equiv 2 \mod ( 3 )$ and $G_1$ is isomorphic to a subgroup of $S_3$ of order multiple of $3$;
\item $G_1$ is contained in a Borel subgroup and it is generated by an element $\sigma$ of order $p$ and an element $g$ of order dividing $2$ such that $\sigma$ and $g$ have the same eigenvector for the eigenvalue $1$.
\end{enumerate}

Moreover, there exist number fields $L_1$, $L_2$, $L_3$ and elliptic curves $\E_1$ defined over $L_1$, $\E_2$ defined over $L_2$, $\E_3$ defined over $L_3$ such that for every $i \in \{ 1, 2, 3 \}$, the $\Gal ( L_i ( \E_i[p] ) / L_i )$-module $\E_i[p]$ is isomorphic to the $G_1$-module $\E[p]$ of the case i. and the local-global divisibility by $p^2$ does not hold over $\E ( L_i )$.
\end{thm}

Clearly case 1. of Theorem \ref{teo1} corresponds to the case when $\E ( k )$ has a point of order $p$ defined over $k$.
Recall that if $ k $ is a number field and $\E$ is an elliptic curve defined over $ k $, then for every $\tau \in \Gal ( k ( \E[p] ) / k )$, the determinant of $\tau$ is the $p$th cyclotomic character.
Then the cases 2. and cases 3. of Theorem \ref{teo1} correspond to the case when $\Q ( \zeta_p + \overline{\zeta_p} ) \subseteq k$. 

By the main result of \cite{PRV1} and Theorem \ref{teo1}, we have the following corollary.

\begin{cor}\label{cor2}
Let $p \geq 5$ be a prime number, let $ k $ be a number field and let $\E$ be an elliptic curve defined over $k$. If $p \equiv 1 \mod ( 3 )$ and $\E$ does not admit any point of order $p$ over $k$, then for every positive integer $n$ the local-global divisibility by $p^n$ holds over $\E ( k )$. If $p \equiv 2 \mod ( 3 )$, $\E$ does not admit any point of order $p$ over $k$ and $[k ( \E[p] ) : k] \neq 3$ and $6$, then for every positive integer $n$ the local-global divisibility by $p^n$ holds over $\E ( k )$.
\end{cor}

\DIM If $k$ does not contain $\Q ( \zeta_p + \overline{\zeta_p} )$, just apply the main result of \cite{PRV1}.
If $p \equiv 1 \mod ( 3 )$ and $k$ contains $\Q ( \zeta_p + \overline{\zeta_p} )$, if there exists $n \in \N$ such that the local-global divisibility by $p^n$ does not hold over $\E ( k )$, then either case 1. or case 3. of Theorem \ref{teo1} holds.
Then $\E$ admits a point of order $p$ defined over $k$.

If $p \equiv 2 \mod ( 3 )$, $\E$ does not admit any point of order $p$ over $k$, and there exists a positive integer $n$ such that the local-global divisibility by $p^n$ does not hold over $\E ( k )$, then case 2. of Theorem \ref{teo1} holds.
Hence $k ( \E[p] ) / k$ is either an extension of degree $3$ or an extension of degree $6$.
\CVD              

\section{Known results}\label{sec11}

In the following proposition we put together the main results of \cite{PRV1} and \cite{PRV2} and we use some results of \cite{GR2}. 

\begin{pro}\label{pro11}
Let $k$ be a number field and let $\E$ be an elliptic curve defined over $k$. Let $p$ be a prime number and, for every $m \in \N$, let $G_m$ be $\Gal ( k ( \E[p^m] ) / k )$. Suppose that there exists $n \in \N$ such that $H^1_{{\rm loc}} ( G_n , \E[p^n] ) \neq 0$. Then one of the following cases holds:
\begin{enumerate}
\item If $p$ does not divide $\vert G_1 \vert$ then either $G_1$ is cyclic of order dividing $p-1$, generated by an element fixing a point of order $p$ of $\E$, or $p \equiv 2 \mod ( 3 )$ and $G_1$ is a group isomorphic either to $S_3$ or to a cyclic group of order $3$;
\item If $p$ divides $\vert G_1 \vert$ then $G_1$ is contained in a Borel subgroup and it is either cyclic of order $p$, or it is generated by an element of order $p$ and an element of order $2$ distinct from $-Id$.
\end{enumerate}
\end{pro}

\DIM Suppose first that $p$ does not divide $\vert G_1 \vert$. 
By the argument in \cite[p. 29]{DZ3}, we have that $G_1$ is isomorphic to its projective image.
By \cite[Proposition 16]{Ser}, then $G_1$ is either cyclic, or dihedral or isomorphic to one of the following groups: $A_4$, $S_4$, $A_5$.

Suppose that the last case holds. 
Then $G_1$ should contain a subgroup isomorphic to $\Z / 2 \Z \times \Z / 2 \Z$ and so $-Id$.
This contradicts the fact that $G_1$ is isomorphic to its projective image.

Suppose that $G_1$ is dihedral. 
Then $G_1$ is generated by $\tau$ and $\sigma$ with $\sigma$ of order $2$ and $\sigma \tau = \tau^{-1} \sigma$.
In particular all the elements of $G_1$ have determinant either $1$ or $-1$.
Suppose that there exists $i \in \N$ such that $\tau^i$ has order dividing $p-1$ and distinct from $1$. 
By \cite[Theorem 2]{GR2}, $\tau^i$ has at least an eigenvalue equal to $1$.
Then, since $\tau^i$ has determinant $-1$, the unique posibility is that $\tau^i$ has order $2$ and so it commutes with $\sigma$.
Thus we get a contradiction because we should have $-Id \in G_1$.
Then $\tau$ has determinant $1$ and order dividing $p+1$.
In particular it has two eigenvalues over $\F_{p^2}$: $\lambda$ and $\lambda^p$.
By \cite[Proposition 17, Lemma 18]{GR2} (or see \cite[Section 3]{CS}) we have that $\tau$ has order $3$.
Then $3$ divides $p+1$ and $G_1$ is isomorphic to $S_3$.

Finally suppose that $G_1$ is cyclic.
If $G_1$ is generated by an element of order dividing $p-1$, by \cite[Theorem 2]{GR2} we have that such an element has an eigenvalue equal to $1$.
On the other hand if the generator of $G_1$ has order not dividing by $p-1$, again by \cite[Proposition 17, Lemma 18]{GR2} we get that such an element has order $3$ and $3$ divides $p+1$. 
\vspace{5 pt}

Suppose now that $p$ divides $\vert G_1 \vert$.
Since $p$ divides the order of $G_1$, by \cite[Proposition 15]{Ser} and the fact that $G_1$ is isomorphic to its projective image, we have that $G_1$ is contained in a Borel subgroup.
In particular the $p$-Sylow subgroup $N$ of $G$ is normal.
Suppose that $G/ N$ is not cyclic.
Then $G_1$ is not isomorphic to its projective image.
Thus $G_1$ is generated by $\sigma$ an element of order $p$, which generates $N$, and $g$ an element of order dividing $p-1$.
Suppose that $g$ has the eigenvalues distinct from $1$.
Then by \cite[Theorem 2]{GR2} (in particular observe that, by \cite[Remark 16]{GR2}, it is not necessary the hypothesis $H^1 ( G_1 , \E[p] ) = 0$), $H^1_{{\rm loc}} ( G_m , \E[p^m] ) = 0$ for every $m \in \N$ and so we get a contradiction.
Then $g$ has an eigenvalue equal to $1$.
Suppose that $g$ has order $\geq 3$.
Then its determinant has order $\geq 3$ and so, since the determinant is the $p$th cyclotomic character, $k$ does not contain $\Q ( \zeta_p + \overline{\zeta_p} )$.  
Then if $g$ and $\sigma$ does not fix the same point of order $p$, by \cite[Theorem 1]{PRV1} we get a contradiction.
On the other hand, since $p$ divides the order of $G_1$, $k ( \E[p] ) \neq k ( \zeta_p )$.
Then by \cite[Theorem 3]{PRV2}, we get a contradiction.

Finally, then either $G_1$ is cyclic generated by an element $\sigma$ of order $p$ or $G_1$ is generated by an element of order $p$ generating a normal subgroup $N$ of $G_1$ and an element $g$ of order $\leq 2$ and which is not $-Id$.
\CVD  

\section{The image of the Galois action over the torsion points}\label{sec12}

In this section we recall some important theorems on the Galois action over the torsion points on an elliptic curve over a number field.
In \cite{GR1} we proved the following lemma, which is a direct consequence of very interesting results of Greicius \cite{Gre} and Zywina \cite{Zyw}.

\begin{lem}\label{lem21}
Given a prime number $p$, a positive integer $n$ and a subgroup $G$ of ${\rm GL}_2 ( \Z/p^n \Z )$, there exists a number field $k$ and an elliptic curve $\E$ defined over $k$ such that there are an isomorphism $\phi \colon \Gal ( k ( \E[p^n] ) / k ) \rightarrow G$ and a $\Z / p^n \Z$-linear homomorphism $\tau \colon \E[p^n] \rightarrow ( \Z / p^n \Z )^2$ such that, for all $\sigma \in \Gal ( k ( \E[p^n] ) / k )$ and $v \in \E[p^n]$, we have $\phi ( \sigma ) \tau ( v ) = \tau ( \sigma ( v ) )$. 
\end{lem}

\DIM See \cite[Lemma 11]{GR1}.
\CVD

The following corollary is an immediate consequence of the previous lemma and \cite[Theorem 3]{DZ3}.

\begin{cor}\label{cor22}
Given a prime number $p$, a positive integer $n$ and a subgroup $G$ of ${\rm GL}_2 ( \Z/p^n \Z )$ such that $H^1_{{\rm loc}} ( G , ( \Z / p^n \Z )^2 ) \neq 0$, there exists a number field $k$ and an elliptic curve $\E$ defined over $k$ such that there exist isomorphisms $\phi \colon \Gal ( k ( \E[p^n] ) / k ) \rightarrow G$ and $\tau \colon \E[p^n] \rightarrow ( \Z / p^n \Z )^2$ such that, for all $\sigma \in \Gal ( k ( \E[p^n] ) / k )$ and $v \in \E[p^n]$, we have $\phi ( \sigma ) \tau ( v ) = \tau ( \sigma ( v ) )$. Then $H^1_{{\rm loc}} ( \Gal ( k ( \E[p^n] ) / k ) , \E[p^n] )$ is isomorphic to $H^1_{{\rm loc}} ( G , ( \Z / p^n \Z )^2 )$ and there exists a finite extension $L$ of $k$ such that $L \cap k ( \E[p^n] ) = k$ and the local-global divisibility by $p^n$ does not hold over $\E ( L )$.
\end{cor}

\section{The prime to $p$ case}\label{sec2}

By Proposition \ref{pro11} we get that if $k$ is a number field and $\E$ is an elliptic curve defined over $k$ such that $\Gal ( k ( \E[p] ) / k )$ has order prime to $p$ and the local-global divisibility by a certain power of $p$ does not hold over $\E ( k )$, then we have substantially two cases to study: the case when $p \equiv 2 \mod ( 3 )$ and $\Gal ( k ( \E[p] ) / k )$ is isomorphic to a subgroup of $S_3$ and the case when $\Gal ( k ( \E[p] ) / k )$ is cyclic of order dividing $p-1$ generated by an element fixing a point of order $p$. 
We do this in the following subsections.

\subsection{The case when $p \equiv 2 \mod ( 3 )$ and the Galois group is a subgroup of $S_3$}

Let $p \equiv 2 \mod ( 3 )$ be a prime number. 
In \cite[Section 5]{GR2} we already found a subgroup $G$ of ${\rm GL}_2 ( \Z / p^2 \Z )$ such that $H^1_{{\rm loc}} ( G, ( \Z / p^2 \Z )^2 ) \neq 0$ and the quotient of $G$ by the subgroup $H$ of the elements congruent to the identity modulo $p$ is a cyclic group of order $3$.  
We use this, the following lemma and the following proposition to extend the example to a group $G^\prime$ containing $G$ such that $G^\prime / H$ is isomorphic to $S_3$.

\begin{lem}\label{lemnuovo}
Let $p$ be a prime number, let $m$ be a positive integer, let $V$ be $( \Z / p^2 \Z )^{2m}$, let $G$ be a subgroup of ${\rm GL}_{2m} ( \Z / p^2 \Z )$ and let $H$ be the subgroup of $G$ of the elements congruent to the identity modulo $p$. 
Then we have the following exact sequence:
\begin{equation}\label{rel31}
0 \rightarrow H^1 ( G/H , V[p] ) \rightarrow H^1 ( G , V[p] ) \rightarrow H^1 ( H, V[p] )^{G/ H} \rightarrow H^2 ( G/H , V[p] ).
\end{equation} 
Moreover the exact sequence 
\[
0 \rightarrow V[p] \rightarrow V \rightarrow V[p] \rightarrow 0
\]
(the first map is the inclusion and the second map the multiplication by $p$) induces the following exact sequence:
\begin{equation}\label{rel32}
H^0 ( G , V[p] ) \rightarrow H^1 ( G , V[p] ) \rightarrow H^1 ( G , V ) \rightarrow H^1 ( G, V[p] ).
\end{equation}
\end{lem}

\DIM Since $H$ is a normal subgroup of $G$, the exact sequence (\ref{rel31}) is just the inflation-restriction sequence for $H$.

Consider the exact sequence of $G$-modules:
\[
0 \rightarrow V[p] \rightarrow V \rightarrow V[p] \rightarrow 0,
\]  
where the first map is the inclusion and the second is the multiplication by $p$.
Then it induces the following exact sequence of $G$-modules
\[
H^0 ( G , V[p] ) \rightarrow H^1 ( G , V[p] ) \rightarrow H^1 ( G , V ) \rightarrow H^1 ( G, V[p] ).
\] 
\CVD

\begin{pro}\label{pro32}
Let $p$ be a prime number, let $m$ be a positive integer, let $V$ be $( \Z / p^2 \Z )^{2m}$ and let $G$ be a subgroup of ${\rm GL}_{2m} ( \Z / p^2 \Z )$ such that:
\begin{enumerate}
\item $G$ has an element $\delta$ not fixing any element of $V$;
\item Let $H$ be the subgroup of $G$ of the element congruent to the identity modulo $p$. $H$ is isomorphic, as $G/H$-module, to a non trivial $G/H$-submodule of $V[p]$;
\item For every $h \in H$ distinct from the identity, the endomorphism $h - Id \colon V / V[p] \rightarrow V / V[p]$ is an isomorphism;
\item $G /H$ has order not divisible by $p$.
\end{enumerate}
Then $H^1_{{\rm loc}} ( G , V ) \neq 0$.
\end{pro}

\DIM Consider the inflation-restriction sequence (see Lemma \ref{lemnuovo})
\[
0 \rightarrow H^1 ( G/H , V[p] ) \rightarrow H^1 ( G , V[p] ) \rightarrow H^1 ( H, V[p] )^{G/ H} \rightarrow H^2 ( G/H , V[p] ).
\]
By Hypothesis 4., $H^1 ( G/H , \A[p] )$ and $H^2 ( G/H , \A[p] )$ are trivial.
Then the restriction $H^1 ( G , V[p] ) \rightarrow H^1 ( H, V[p] )^{G/ H}$ is an isomorphism.
Since the action of $H$ over $V[p]$ is trivial and $H$ is an abelian group of exponent $p$, we have that $H^1 ( H , V[p] )^{G/H}$ is isomorphic to ${\rm Hom}_{\Z / p \Z[G/ H]} ( H, V[p] )$. 
By hypothesis 2., $H$ is isomorphic to a non trivial $G/H$-submodule of $V[p]$.
Then there exists $\phi \colon H \rightarrow V[p]$ an injective homomorphism of $\Z/ p \Z[G/ H]$-modules.
Let $[Z]$ be in $H^1 ( G , V[p] )$ such that its image in $H^1 ( H, V[p] )^{G/ H}$ is the class of $\phi$.
In particular observe that since $\phi$ is an injective homomorphism, $[Z] \neq 0$.

Now observe that $H^0 ( G , V[p] ) = 0$ by hypothesis 1..
Then, by Lemma \ref{lemnuovo}, we have the following exact sequence of $G$-modules
\[
0 \rightarrow H^1 ( G , V[p] ) \rightarrow H^1 ( G , V ) \rightarrow H^1 ( G, V[p] ).
\] 
Let us call $[W] \in H^1 ( G , V )$ the image of $[Z] \in H^1 ( G , V[p] )$ defined above by the injective map $H^1 ( G , V[p] ) \rightarrow H^1 ( G , V )$.
Since $[Z] \neq 0$, the same holds for $[W]$.
Moreover, since $G/ H$ is not divisible by $p$, the restriction $H^1 ( G , V ) \rightarrow H^1 ( H , V )$ is injective.
Since, by 3., the image by the restriction of $[W]$ over $H^1 ( H , V )$ is in $H^1_{{\rm loc}} ( H , V )$, we have that $[W]$ is a non-trivial element of $H^1_{{\rm loc}} ( G , V )$.
\CVD  

\begin{cor}\label{cor33}
Let $p$ be an odd prime such that $p \equiv 2 \mod ( 3 )$. Let $G$ be the subgroup of ${\rm GL}_2 ( \Z / p^2 \Z )$, such that $G$ is generated by
\[
\tau = 
\left( 
\begin{array}{cc} 
1 & -3 \\
1 & -2 \\
\end{array}
\right)
\]
of order $3$, $\sigma$ of order $2$ such that $\sigma \tau \sigma^{-1} = \tau^2$ and
\[
H = \bigg\{
\left( 
\begin{array}{cc} 
1 + p ( a - 2b ) & 3p ( b-a ) \\  
-pb & 1 - p ( a - 2b ) \\
\end{array}
\right), \ a, b \in \Z/ p^2 \Z
\bigg\}.
\]
Then $H^1_{{\rm loc}} ( G , ( \Z / p^2 \Z )^2 ) \neq 0$.
\end{cor}

\DIM Observe that conditions 1. and 4. of Proposition\ref{pro32} hold for $G$.
Moreover condition 3. holds by \cite[Section 5]{GR2}.
Observe that $G/ H$ is isomorphic to $S_3$ and recall that $S_3$ has a unique irreducible representation of dimension $2$ over $\F_p$.
Then condition 2. of Proposition \ref{pro32} is equivalent to prove that $H$ is stable by the conjugation by $\tau$ and $\sigma$.
In \cite[Section 5]{GR2} we proved that the conjugation by $\tau$ sends $H$ to $H$.
A straightforward computation shows that if $\overline{\sigma}$ has order $2$ in $G/ H$ and $\overline{\sigma} \overline{\tau} \overline{\sigma}^{-1} = \overline{\tau}^2$, then there exists $\alpha , \beta \in \F_p$ such that   
\[
\overline{\sigma} =
\left( 
\begin{array}{cc} 
\alpha - 2\beta  & 3 ( \beta-\alpha ) \\  
\beta & 2 \beta - \alpha \\
\end{array}
\right).
\] 
Another straightforward computation shows that $\sigma H \sigma^{-1} = H$ (since $H$ has dimension $2$ as $\F_p$-vector space, it is sufficient to verify this over a basis and it is a trivial  verification at least over the element of $H$ with $\alpha = a$ and $\beta = b$).
Then, by Proposition \ref{pro32}, we have $H^1_{{\rm loc}} ( G , ( \Z / p^2 \Z )^2 ) \neq 0$.
\CVD

\subsection{The case when the Galois group is cyclic of order dividing $p-1$}

By using some results of the previous subsection, we study the second case. 

\begin{lem}\label{lem34}
Let $p$ be a prime number and let $V$ be $( \Z/ p^2 \Z )^2$.
Let $\lambda \in ( \Z/ p^2 \Z )^\ast$ of order dividing $p-1$ and let $G$ be the following subgroup of ${\rm GL}_2 ( V )$:
\[
G =
\bigg\langle
g = 
\left( 
\begin{array}{cc} 
\lambda & 0 \\
0 & 1 \\
\end{array}
\right),
h ( 1, 0 ) =
\left( 
\begin{array}{cc} 
1 + p & 0 \\
0 & 1 - p \\
\end{array}
\right),
h ( 0, 1 ) =
\left( 
\begin{array}{cc} 
1 & p \\
0 & 1 \\
\end{array}
\right)    
\bigg\rangle.
\]
Then $H^1_{{\rm loc}} ( G, V ) \neq 0$.
\end{lem}

\DIM Observe that the subgroup $H$ of $G$ of the elements congruent to the identity modulo $p$ is the group generated by $h ( 1, 0 )$ and $h ( 0, 1 )$.
Since $G / H$ has order not divisible by $p$, $H^1 ( G/ H , V[p] ) = 0$ and $H^2 ( G/ H , V[p] ) = 0$.
Then, from the exact sequence (\ref{rel31}) in Lemma \ref{lemnuovo}, we get an isomorphism from $H^1 ( G/ H , V[p] )$ to $H^1 ( H , V[p] )^{G/ H}$.
Since $H$ acts like the identity over $V[p]$ and $V[p]$ and $H$ are abelian groups with exponent $p$, $H^1 ( H , V[p] )^{G/ H} = {\rm Hom}_{\Z / p \Z[G/H]} ( H , V[p] )$.
Observe that $g h ( 0, 1 ) g^{-1} = h ( 0, 1 )^\lambda$ and $g ( p , 0 ) = \lambda ( p, 0 )$.
Then we can define a non-trivial $\Z / p \Z[G/H]$ homomorphism $\phi$ from $H$ to $V[p]$ by sending $h ( 0, 1 )$ to $( p, 0 )$ and $h ( 1, 0 )$ to $( 0, 0 )$ and extending it by linearity. 
Let $Z$ be a cocycle representing the class $[Z]$ in $H^1 ( G , V[p] )$ corresponding to $\phi$. 
By the exact sequence (\ref{rel32}) in Lemma \ref{lemnuovo} we get an homomorphism from $H^1 ( G, V[p] )$ to $H^1 ( G , V )$.
Let $W$ be a cocycle representing the class $[W]$ in $H^1 ( G , V )$ image of $[Z]$ for such homomorphism.
Let us show that $[W] \in H^1_{{\rm loc}} ( G, V )$ and $[W] \neq 0$.
Since $G / H$ has order not divisible by $p$, it is sufficient to prove that the image by the restriction of $[W]$ to $H^1 ( H , V )$ is in $H^1_{{\rm loc}} ( H , V )$. 
For all $a , b$ integers put $h ( a , b ) = a h ( 1, 0 ) + b h ( 0, 1 )$.
Then, by definition of $[Z]$, $h ( a, b )$ is sent to $( b p , 0 )$.
An easy calculation shows that for every $a , b$, there exists $x$, $y$ in $\Z / p^2 \Z$ such that $( h - Id ) ( x, y ) = ( bp, 0 )$.
This proves that $[W] \in H^1_{{\rm loc}} ( G, V )$.

Finally observe that for every $x$, $y$ in $\Z / p^2 \Z$ such that $( h ( 1, 0 ) - Id ) ( x, y ) = ( 0, 0 )$, we have $x \equiv 0 \mod ( p )$ and $y \equiv 0 \mod ( p )$.
On the other hand for every $x$, $y$ in $\Z / p^2 \Z$ such that $( h ( 1, 0 ) - Id ) ( x, y ) = ( p, 0 )$, we have $y \equiv 1 \mod ( p )$.
Thus $[W] \neq 0$.
\CVD

\section{The case when $p$ divides the order of the Galois group}\label{sec4}

By Proposition \ref{pro11} we get that if $k$ is a number field and $\E$ is an elliptic curve defined over $k$ such that $\Gal ( k ( \E[p] ) / k )$ has order divisible by $p$ and the local-global divisibility by a certain power of $p$ does not hold over $\E ( k )$, then $\Gal ( k ( \E[p] ) / k )$ is contained in a Borel subgroup and it is generated by an element $\sigma$ of order $p$ and an element $g$ of order dividing $2$ and distinct from $-Id$.
In the following subsections we study first the case when $g$ and $\sigma$ fix the same element of order $p$, then the case when $g$ and $\sigma$ do not fix any element of order $p$.

\subsection{The case when $g$ and $\sigma$ fix the same vector}

In this section we prove the following result.

\begin{lem}\label{lem42}
Let $V$ be $( \Z / p^2 \Z )^2$ and let $G$ be the following subgroup of ${\rm GL}_2 ( \Z / p^2 \Z )$:
\[
G =
\bigg\langle
g = 
\left( 
\begin{array}{cc} 
1 & 0 \\
0 & -1 \\
\end{array}
\right),
\sigma =
\left( 
\begin{array}{cc} 
1 + p & 1 \\
2  p & 1 + p \\
\end{array}
\right),
h =
\left( 
\begin{array}{cc} 
1 + p & 0 \\
0 & 1 - p \\
\end{array}
\right)    
\bigg\rangle.
\]
Then $H^1_{{\rm loc}} ( G, V ) \neq 0$.
\end{lem}

\DIM Let $H$ be the subgroup of $G$ of the elements congruent to $1$ modulo $p$.
Let $\overline{g}$ and $\overline{\sigma}$ be the classes of $g$ and $\sigma$ modulo $H$.
We have that $H^1 ( G /H , V[p] ) \neq 0$.
In fact we can define a cocycle $Z \colon G/ H \rightarrow V[p]$, which is not a coboundary, by sending for every integer $i_1 , i_2$, $Z_{\overline{g}^{i_1} \overline{\sigma}^{i_2}} = ( p i_2 ( i_2 - 1 ) / 2 , ( -1 )^{i_1} p i_2 )$.
Since $H$ is normal, we have an injective homorphism (the inflation) from $H^1 ( G/ H , V[p] )$ to $H^1 ( G , V[p] )$.
By abuse of notation we still call $Z$ a cocycle representing the image of the class of $Z$ in $H^1 ( G , V[p] )$.
Moreover, see Lemma \ref{lemnuovo} in particular the sequence (\ref{rel32}), we have a homomorphism from $ 
H^1 ( G , V[p] )$ to $H^1 ( G , V )$.
It sends the class of $Z$ in $H^1 ( G, V[p] )$ in the class of the cocyle $W$ representing the class $[W] \in H^1 ( G , V )$.
We shall prove that $[W] \in H^1_{{\rm loc}} ( G, V )$ and $[W] \neq 0$.

First of all let us observe that for every $a, b, c, d \in \Z / p^2 \Z$, we have
\[
\left( 
\begin{array}{cc} 
1 + ap & 1 + bp \\
cp & 1 + dp \\
\end{array}
\right)^p =
\left( 
\begin{array}{cc} 
1 & p \\
0 & 1  \\
\end{array}
\right).
\]
To verify this write 
\[
\left( 
\begin{array}{cc} 
1 + ap & 1 + bp \\
cp & 1 + dp \\
\end{array}
\right) = 
\left( 
\begin{array}{cc} 
1 & 0 \\
0 & 1  \\
\end{array}
\right) +
\left( 
\begin{array}{cc} 
ap & 1 + bp \\
cp & dp \\
\end{array}
\right)
\]
and observe that
\[
\left( 
\begin{array}{cc} 
ap & 1 + bp \\
cp & dp \\
\end{array}
\right)^2 \equiv
\left( 
\begin{array}{cc} 
0 & 0 \\
0 & 0 \\
\end{array}
\right) \mod ( p ), \
\left( 
\begin{array}{cc} 
ap & 1 + bp \\
cp & dp \\
\end{array}
\right)^4 =
\left( 
\begin{array}{cc} 
0 & 0 \\
0 & 0 \\
\end{array}
\right).     
\]
Thus the subgroup $H$ of $G$ of the elements congruent to the identity modulo $p$ is 
\[
H = \bigg \langle
\left(  
\begin{array}{cc} 
1 & p \\
0 & 1  \\
\end{array}
\right),
\left( 
\begin{array}{cc} 
1 + p & 0 \\
0 & 1 - p \\
\end{array}
\right)
\bigg \rangle.    
\]
Now observe that, since $H$ and $\langle \sigma , H \rangle$ are normal in $G$, for every $\tau \in G$ there exist integers $i_1 , i_2 , i_3$ and $h \in H$ such that $\tau = g^{i_1} \sigma^{i_2} h^{i_3}$.
By definition of $W$, $W_\tau = ( p ( i_2 - 1 ) , ( -1 )^{i_1} p i_2 )$.
If $i_2 \equiv 0 \mod ( p )$, then clearly $W_\tau = ( 0, 0 )$ and so $W_\tau = ( \tau - Id ) ( ( 0, 0 ) )$.
Then we can suppose $i_2 \not \equiv 0 \mod ( p )$.
It is simple to prove by induction on $i_2$ that
\[
\sigma^{i_2} =
\left( 
\begin{array}{cc} 
1 + ap & i_2 + b p \\
2 i_2 p & 1 + c p \\
\end{array}
\right)  
\]
for certain $a, b, c \in \Z / p^2 \Z$.
Moreover  $\sigma^{i_2} h^{i_3}$ has still the coefficient at the top on the right congruent to $i_2$ modulo $p$ and the low coefficient on the left equal to $2 i_2 p$.
From these remarks is an easy exercise to prove that there exist $\alpha$ and $\beta \in \Z / p^2 \Z$ such that $W_\tau = ( \tau - Id ) ( ( \alpha , p \beta ) )$.
Then $[W]$ is in $H^1_{{\rm loc}} ( G, V )$.

Finally let us observe that $W$ is not a coboundary.
Let $\alpha$, $\beta \in \Z / p^2 \Z$ such that $W_\sigma = ( 0, p ) = ( \sigma - Id ) ( ( \alpha , \beta ) )$.
Then $\alpha \not \equiv 0 \mod ( p )$.
On the other hand let $h$ be in $H$ be  
\[
h =    
\left( 
\begin{array}{cc} 
1 + p & 0 \\
0 & 1 - p \\
\end{array}
\right).
\]  
Then $W_h = ( 0, 0 )$ and so for every $\alpha , \beta \in \Z / p^2 \Z$ such that $( h - id ) ( ( \alpha , \beta ) ) = ( 0, 0 )$, we have $\alpha \equiv 0 \mod ( p )$.
Hence $W$ is not a coboundary.
\CVD

\begin{obs}\label{remnuovina} 
Since we shall use it in the next subsection, we make the following remark.
In the proof of the previous lemma, we observed that for every $a, b, c, d \in \Z / p^2 \Z$, we have
\[
\left( 
\begin{array}{cc} 
1 + ap & 1 + bp \\
cp & 1 + dp \\
\end{array}
\right)^p =
\left( 
\begin{array}{cc} 
1 & p \\
0 & 1  \\
\end{array}
\right).
\]
In a similar way, for every integer $m \geq 2$, and every $a_m, b_m, c_m, d_m \in \Z/ p^m \Z$, we have
\[
\left( 
\begin{array}{cc} 
1 + a_mp & 1 + b_mp \\
c_mp & 1 + d_mp \\
\end{array}
\right)^{p^{m-1}} =
\left( 
\begin{array}{cc} 
1 & p^m \\
0 & 1  \\
\end{array}
\right).
\]
\end{obs} 

\begin{cor}\label{cornuovo1}
Let $V$ be $( \Z / p^2 \Z )^2$ and let $G$ be the following subgroup of ${\rm GL}_2 ( \Z / p^2 \Z )$:
\[
\widetilde{G} =
\bigg\langle
\sigma =
\left( 
\begin{array}{cc} 
1 + p & 1 \\
2  p & 1 + p \\
\end{array}
\right),
h =
\left( 
\begin{array}{cc} 
1 + p & 0 \\
0 & 1 - p \\
\end{array}
\right)    
\bigg\rangle.
\]
Then $H^1_{{\rm loc}} ( \widetilde{G}, V ) \neq 0$
\end{cor}

\DIM Observe that $\widetilde{G}$ is a subgroup of index $2$ of the group $G$ of Lemma \ref{lem42}.
Since $p \neq 2$ the restriction $H^1_{{\rm loc}} ( G, V ) \rightarrow H^1_{{\rm loc}} ( \widetilde{G}, V )$ is injective.
Thus, since $H^1_{{\rm loc}} ( G, V ) \neq 0$, $H^1_{{\rm loc}} ( \widetilde{G}, V ) \neq 0$.
\CVD

\subsection{The case when $g$ and $\sigma$ does not fix the same vector and proof of Theorem \ref{teo1}}

In this section we study the case when $g$ and $\sigma$ have a common eigenvector, but they do not fix the same vector. 
 
\begin{lem}\label{lem43}
Let $n \in \N$ and let $G$ be a subgroup of ${\rm GL}_2 ( \Z / p^n \Z )$. Let $H$ be the subgroup of $G$ of the elements congruent to the identity modulo $p$. Suppose that $G / H$ is contained in a Borel subgroup and it is generated by an element $g$ of order $2$ and an element $\sigma$ of order $p$ such that $\sigma$ and $g$ do not fix the same element of order $p$. Let $V_n$ be $( \Z / p^n \Z )^2$. Then $H^1_{{\rm loc}} ( G, V_n ) = 0$.
\end{lem}

\DIM  By repalcing $V$ with $V_n$, by observing that $H^0 ( G , V_n[p^{n-1}] ) = 0$ because the group generated by $g$ and $\sigma$ do not fix any element of $V_n[p^{n-1}]$, and by copying the proof of Lemma \ref{lemnuovo}, we get the following exact sequence
\begin{equation}\label{relsa}
0 \rightarrow H^1 ( G, V_n[p] ) \rightarrow H^1 ( G, V_n ) \rightarrow H^1 ( G, V_n[p^{n-1}] ).
\end{equation}
Suppose that $H^1_{{\rm loc}} ( G, V_n ) \neq 0$.
Then $H^1_{{\rm loc}} ( G, V_n )[p] \neq 0$ and let $Z$ be a cocycle representing a non-trivial class $[Z] \in H^1_{{\rm loc}} ( G, V_n )[p]$.
Let us observe that $[Z]$ is in the kernel of $H^1 ( G, V_n ) \rightarrow H^1 ( G, V_n[p^{n-1}] )$ (here we generalize the proof of \cite[Lemma 13]{GR2}).
Since $[Z]$ has order $p$, then $pZ$ is a coboundary and so there exists $v \in V_n$ such that, for every $\tau \in G$, $pZ_\tau = \tau ( v ) - v$.
Let us observe that $v \in V_n[p^{n-1}]$.
Since for every $\tau$, $\tau ( v ) - v \in V_n[p^{n-1}]$, we have $v \in \cap_{\tau \in G} \ker ( p^{n-1} ( \tau - Id ) )$.
Since $G$ does not fix any element of order $p$ the unique possibility is that $v \in V_n[p^{n-1}]$.
Then (see the sequence (\ref{relsa})) $[Z]$ is in the image of $H^1 ( G, V_n[p] ) \rightarrow H^1 ( G, V_n )$.
By abuse of notation we call $[Z]$ the class in $H^1 ( G, V_n[p] )$ sent to $[Z]$.

Consider now the inflation-restriction sequence
\begin{equation}\label{relsa2}
0 \rightarrow H^1 ( G / H , V_n[p] ) \rightarrow H^1 ( G, V_n[p] ) \rightarrow H^1 ( H , V_n[p] )^{G/ H}.
\end{equation}
Let us observe that $H^1 ( G/ H , V_n[p] ) = 0$.
Let $W \colon G/ H \rightarrow V_n[p]$ be a cocycle.
Since $\sigma$ and $g$ are contained in a Borel subgroup, $g$ has order $2$, and $g$ and $\sigma$ do not fix any non-zero element of $V_n[p^{n-1}]$, we can choose a basis of $V_n$ such that $( p^{n-1} , 0 )$ is fixed by $\sigma$ and $g ( ( p^{n-1} , 0 ) ) = ( - p^{n-1} , 0 )$ and $( 0 , p^{n-1} )$ is sent to $( p^{n-1} , p^{n-1} )$ by $\sigma$ and fixed by $g$.
Observe that, since summing a coboundary to $W$ does not change its class, we can suppose that $W_\sigma = ( 0 , p^{n-1} )$.
Then, for every integer $i$, $W_{\sigma^i} = ( p^{n-1} i ( i-1 ) / 2 , p^{n-1} i )$.
Observe that since $g$ has order $2$, $W_{g^2} = W_g + g W_g = ( 0, 0 )$.
In particular there exists $a \in \Z/ p^n \Z$ such that $W_g = ( p^{n-1} a , 0 )$, which is fixed by $\sigma$.
Thus $W_{g \sigma g^{-1}} = g W_\sigma = ( p^{n-1} , -p^{n-1} )$.
On the other hand $g \sigma g^{-1} = \sigma^{-1}$ and so $W_{\sigma^{-1}} = ( - p^{n-1} , -p^{n-1} )$.
We then get a contradiction.
Thus, by the sequence (\ref{relsa2}), we get that to every class of $H^1 ( G, V_n[p] )$ we can associate a class in $H^1 ( H , V_n[p] )^{G/ H}$.
Since $H$ acts like the identity over $V_n[p]$, we have that $H^1 ( H , V_n[p] )^{G/ H}$ is a subgroup of ${\rm Hom} ( H , V_n[p] )$.
In particular we can associate to $[Z] \in H^1 ( G , V[p] )$ defined above a homomorphism from $H$ to $V_n[p]$.
We now need the following lemma.

\begin{lem}\label{lem44}
Let $\tau$ be in $H$ and let $\sigma_n \in G$ be such that $\sigma_n$ is sent to $\sigma$ by the projection of $G$ over $G/ H$. 
Then there exists $\tau_d, \tau_l \in H$, $\lambda \in \N$, such that $\tau_d$ is diagonal, $\tau_l$ is lower unitriangular and $\tau = \tau_d \tau_l \sigma_n^{p \lambda}$.
Then $H$ is generated by its subgroups of the diagonal matrices, its subgroup of the lower unitriangular matrices and $\sigma_n^p$.
\end{lem} 

\DIM We remark that $\sigma_n^p \in H$.
In fact $\sigma_n^p \equiv Id \mod ( p )$. 

We first show that every $\tau \in H$ can be written as a product of a lower triangular matrix $\tau_L \in H$ and a power of $\sigma_n^p$.
Since $\tau \in H$, $\tau \equiv Id \mod ( p )$ and so there exist $e, g, m, r \in \Z / p^n \Z$ such that 
\[
\tau = \left(
\begin{array}{cc}
1 + pe & pg \\
pm & 1 + pr\\
\end{array}
\right).
\]
We prove by induction that for every integer $i \geq 1$, there exists $\lambda_i \in \Z / p^n \Z$ such that 
\begin{equation}\label{rel23}
\tau \sigma_n^{p \lambda_i} =
\left(
\begin{array}{cc}
1 + pe_i & p^i g_i \\
pm_i & 1 + pr_i\\
\end{array}
\right)
\end{equation}
for certain $e_i, g_i, m_i, r_i \in \Z /p^n \Z$.
If $i = 1$ then for $\lambda_1 = 0$ the relation (\ref{rel23}) is satisfied.
Suppose that (\ref{rel23}) is satisfied for an integer $i \geq 1$.
Then there exists $\lambda_i \in \Z / p^n \Z$ such that
\[
\tau \sigma_n^{p \lambda_i} =
\left(
\begin{array}{cc}
1 + pe_i & p^i g_i \\
pm_i & 1 + pr_i\\
\end{array}
\right)
\]
for certain $e_i, g_i, m_i, r_i \in \Z /p^n \Z$.
Set $\lambda_{i+1}$ an element of $\Z /p^n \Z$ such that $p \lambda_{i+1} = p \lambda_i - p^i g_i$.
Observe that this element exists because $i \geq 1$.
By Remark \ref{remnuovina} we have
\begin{align*}
\sigma_n^{-p^i g_i} & = 
\left(
\begin{array}{cc}
1 + p^{i+1}a_{i+1} & p^i + p^{i+1}b_{i+1} \\
p^{i+1}c_{i+1} & 1 + p^{i+1}d_{i+1}\\
\end{array}
\right)^{-g_i} \\
& = \left(
\begin{array}{cc}
1 + p^{i+1}a^\prime_{i+1} & -p^i g_i + p^{i+1}b^\prime_{i+1} \\
p^{i+1}c^\prime_{i+1} & 1 + p^{i+1}d^\prime_{i+1}\\
\end{array}
\right),
\end{align*}
for certain $a^\prime_{i+1}, b^\prime_{i+1}, c^\prime_{i+1}, d^\prime_{i+1} \in \Z /p^n \Z$.
By a short computation
\begin{align*}
\tau \sigma_n^{p \lambda_{i+1}} & = \tau \sigma_n^{p\lambda_i}\sigma_n^{-p^i g_i} \\
& =
\left(
\begin{array}{cc}
1 + pe_i & p^i + p^ig_i \\
pm_i & 1 + pr_i\\
\end{array}
\right) 
\left(
\begin{array}{cc}
1 + p^{i+1}a^\prime_{i+1} & -p^i g_i + p^{i+1}b^\prime_{i+1} \\
p^{i+1}c^\prime_{i+1} & 1 + p^{i+1}d^\prime_{i+1}\\
\end{array}
\right) \\
& = \left(
\begin{array}{cc}
1 + pe_{i+1} &  + p^{i+1}g_{i+1} \\
pm_{i+1} & 1 + pr_{i+1}\\
\end{array}
\right),
\end{align*}
for certain $e_{i+1}, g_{i+1}, m_{i+1}, r_{i+1} \in \Z /p^n \Z$.
Then (\ref{rel23}) is verified for $\lambda_{i+1}$ that satisfies $p\lambda_{i+1} = p\lambda_i -p^i g_i$.
In particular for $i = n$ we have
\[
\tau \sigma_n^{p \lambda_n} =
\left(
\begin{array}{cc}
1 + pe_n & 0 \\
pm_n & 1 + pr_n\\
\end{array}
\right).
\]
Then seting $\tau_L = \tau \sigma_n^{p \lambda_n}$ and $\lambda = -\lambda_n$, we have shown that $\tau$ can be written as a product of a lower triangular matrix $\tau_L \in H$ and the power $\sigma_n^{p \lambda}$ of $\sigma_n^p$.

Observe that, to conclude the proof, it is sufficient to show that $\tau_L$ can be written as the product of a diagonal matrix $\tau_d \in H$ and a lower unitriangular matrix $\tau_l \in H$.
Let $g_n \in G$ be an element of order $2$ whose projection to $G/ H$ is $g$. 
Since $H$ is normal in $G$, $g_n \tau_L g_n^{-1} \in H$.
Then $g_n \tau_L g_n^{-1} \tau_L^{-1} \in H$.
Moreover by a simple computation, we have
\[
g_n \tau_L g_n^{-1} \tau_L^{-1} =
\left(
\begin{array}{cc}
1 & 0 \\
-2pm_n / ( pe_n + 1 ) & 1 \\
\end{array}
\right).
\]
Thus
\[
( g_n \tau_L g_n^{-1} \tau_L^{-1} )^{- ( pe_n + 1 ) / 2 ( pr_n + 1 )} =
\left(
\begin{array}{cc}
1 & 0 \\
pm_n / ( pr_n + 1 ) & 1 \\
\end{array}
\right) \in H. 
\]
Call such a matrix $\tau_l$ and observe that
\[
\tau_L =
\left(
\begin{array}{cc}
1 + pe_n & 0 \\
0 & 1 + pr_n\\
\end{array}
\right)
\left(
\begin{array}{cc}
1 & 0 \\
pm_n / ( pr_n + 1 ) & 1 \\
\end{array}
\right).
\]
Call the diagonal matrix $\tau_d$.
Since $\tau_L, \tau_l \in H$, also $\tau_d \in H$, proving the assumption.
\CVD

Let $\tau \in H$.
By Lemma \ref{lem44} there exists $\tau_l \in H$ a lower unitriangular matrix, $\tau_D \in H$ a diagonal matrix and $\lambda \in \Z$ such that $\tau = \tau_l \tau_D \sigma_n^{\lambda p}$ and consider the homorphism associated to $[Z] \in H^1 ( G, V_n[p] )$.
Since the cocycle $Z$ have values in $V_n[p]$, in particular $Z_{\sigma_n} \in V_n[p]$ and, by properties of cocycle, $Z_{\sigma_n^p} = ( 0, 0 )$.
On the other hand, since $g_n \tau_D g_n^{-1} = \tau_D$, we have that there exists $b \in \Z / p^n \Z$ such that $Z_{\tau_D} = ( 0, p^{n-1} b )$.
Since if $p^{n-1} b$ is distinct from $0$, $( 0, p^{n-1} b )$ generates $V[p]$ as $G/ H$-module and $g_n \tau_l g_n^{-1} = \tau_l^{-1}$, there exists $a \in \Z/ p^n \Z$ such that $Z_{\tau_l} = ( p^{n-1} a , 0 )$.
Observe that for every $( \alpha , \beta ) \in V_n$, $( \tau_l - Id ) ( \alpha , \beta ) = ( p^{n-1} a , 0 )$ only if $p^{n-1} a = 0$.
Then if the image of $Z$ satisfies the local conditions over $V_n$, the homomorphism associated to $Z$ is trivial and so $Z$ is a coboundary.
\CVD

{\bf Proof of Theorem \ref{teo1}}. Let $k$ be a number field and let $\E$ be an elliptic curve defined over $k$.
Set, as before, for every $m \in \N$, $G_m$ the group $\Gal ( k ( \E[p^m] ) / k )$.
By Proposition \ref{pro11} and Lemma \ref{lem43} either case 1., or case 2., or case 3. holds. 

By Corollary \ref{cor22} and Lemma \ref{lem34} we can find a number field $L_1$ and an elliptic curve $\E_1$ defined over $L_1$ such that $\Gal ( k ( \E_1[p] ) / k )$ is cyclic of order dividing $p-1$ and its generator has an eigenvalue equal to $1$

Suppose that $p \equiv 2 \mod ( 3 )$.
By Corollary \ref{cor22} and Corollary \ref{cor33} we can find a number field $L_2$ and an elliptic curve $\E_2$ defined over $L_2$ such that $\Gal ( k ( \E_2[p] ) / k )$ is isomorphic to a subgroup of $S_3$ of order multiple of $3$ and the local-global divisibility by $p^2$ does not hold over $\E ( L_2 )$.

Finally, by Corollary \ref{cor22}, Lemma \ref{lem42} and Corollary \ref{cornuovo1}, we can find a number field $L_3$ and an elliptic curve $\E_3$ defined over $L_3$ such that $\Gal ( k ( \E_3[p] ) / k )$ is contained in a Borel subgroup and generated by $\sigma$ of order $p$ and $g$ of order $2$, fixing the same point of $\E_3[p]$ of order $p$.
\CVD

\end{document}